\documentclass[12pt]{amsart}
\usepackage{amsmath,amscd,latexsym,verbatim,amssymb}
\usepackage{times}

\newcommand{\cf}{{\it cf.} }
\newcommand{\ie}{{\it i.e.} }
\newcommand{\eg}{{\it e.g.} }

\newcommand{\Q}{\mathbf{Q}}

\newcommand{\C}{\mathbf{C}}

  \newcommand{\sA}{{\mathcal{A}}}

\newcommand{\sM}{{\mathcal{M}}}

\renewcommand{\epsilon}{\varepsilon}
\renewcommand{\phi}{\varphi}

\renewcommand{\lim}{\varprojlim}

\font\smit= cmti10 at 10pt

\newcounter{spec}
{\end{list}}

\swapnumbers

\theoremstyle{definition}

\numberwithin{equation}{section}

 at10pt
 at12pt

\begin{document}

\title{G-functions, motives, and unlikely intersections - old and new.}
\author{Yves
Andr\'e}
  \address{Institut Math\'ematique de Jussieu, Sorbonne Universit\'e,
  Paris 05\\France.}
\email{yves.andre@imj-prg.fr}
  \date{20/08/2024}
\maketitle
% \centerline{\smit dedicated to Enrico Bombieri, with admiration, for his 85th birthday}

  \begin{abstract}  In this survey, we outline the role of G-functions in arithmetic geometry, notably their link with Picard-Fuchs differential equations and periods. We explain how polynomial relations between special values of G-functions arising from a pencil of algebraic varieties may occur at a parameter where the fiber has more ``motivic" symmetries; and how Bombieri's principle of global relations can be used to control the height of such parameters (which was also one of the origins of the Andr\'e-Oort conjecture). At the end, we sketch the recent revival of the G-function method in the context of unlikely intersections and the Zilber-Pink conjecture.   
  \end{abstract}

%\tableofcontents

%%\renewcommand{\abstractname}{Summary}

\begin{sloppypar} 
  
  \section{G-functions and G-operators.}(General references: \cite{A1}\cite{DGS}\cite{R}])\label{s1} 
\subsection{G-functions} G-functions were introduced by C.-L. Siegel in his great memoir on diophantine approximation \cite{S}. 
They are formal power series 
$$ g(z) = \sum_0^\infty \,a_n z^n$$
with coefficients in some number field $K$, which satisfy a linear differential equation with polynomial coefficients, and such that the conjugates of the $a_n$'s as well as their denominators have at most an exponential growth in $n$. 

At the end of his memoir, Siegel stated without proof an irrationality property of some special values of G-functions. It took more than 50 years to get a proof of Siegel's claim in full generality: this was achieved in E. Bombieri's pioneering paper ``On G-functions" \cite{B}. 
 
 \subsection{G-operators and geometry.}\label{s2}
 Bombieri's paper initiates the general study of linear differential operators satisfied by G-functions, building on $p$-adic techniques. Specifically, it focuses on differential operators $\Lambda \in K\langle z, d/dz\rangle$ satisfying Bombieri's condition: the product $\prod_p R_p(\Lambda)$ of the $p$-adic radii of convergence at the generic disk is non-zero. 
 
 Writing $\Lambda$ in matrix form $\frac{d}{dz} - A(z)$, let $T(z)\in K[z]$ be a common denominator for the entries of $A$, and $\frac{1}{n!}(\frac{d }{d z  })^n  - A_n(z)$ stand for the unique left multiple of $\Lambda$ of this form. It turns out that Bombieri's condition is equivalent to the condition (introduced by A. Galochkin) that the common denominator of the coefficients of the entries of $T(z)^n A_n(z)$ has an exponential growth in $n$. Such operators are called {\it G-operators}.  
  They are fuchsian, with rational exponents (\cf \cite{K}).
 
 An important theorem due to Chudnovsky asserts that {\it for any G-function $g(z)$, a differential operator of minimal order which kills $g(z)$ is always a G-operator}\footnote{Bombieri's proof of Siegel irrationality claim is conditional to the Bombieri's condition, but Chudnovsky's theorem allows to take it for granted}. Conversely, any formal power series solution of a G-operator is a G-function\footnote{if $0$ is a singularity, solutions involve powers of $\log$ and fractional powers of $z$, but a similar statement also holds for the power series which ``occur" in the solutions. }. 
 
 \smallskip The significance of G-operators lies in their conjectural geometric origin (Bombieri-Dwork). Here, ``$\Lambda $ is of geometric origin" can be understood as follows: $\Lambda  \in K\langle z, d/dz\rangle$ is a product of factors of Picard-Fuchs operators; equivalently, the associated differential module is an iterated extension of factors of differential modules of Gauss-Manin type. 
 
 It was proven in \cite[V App.]{A1} that {\it differential operators of geometric origin are G-operators}. The converse remains very much open. 
 
 The theory generalizes to several variables \cite{AB}\cite{D}. 
 
 \section{G-values and periods.}(General references: \cite{A1}\cite{A6}\cite{F}\cite{FN})\label{s2} 

  \subsection{Periods}  Integrals of regular $1$-forms on an algebraic curve depend on the integration path, and the ambiguities become periods of abelian functions through Abel-Jacobi inversion. In spite of the lack of such an interpretation in higher dimensions, the term ``period" happened to be used, more generally, for any integral $\,\int_\Delta \omega\,$ of an algebraic $n$-form $\omega$ taken on a domain $\Delta$ bounded by algebraic equations.
  
  In a tradition which goes back to Euler, Legendre and Gauss, two special cases are prominent: 
  
  - the arithmetic case (formalized in \cite{KZ}): (arithmetic) periods are complex numbers whose real and imaginary parts are integrals $\int_\Delta \omega$ of an algebraic differential form $\omega$ on an algebraic variety $X$ defined over $\mathbb Q$, taken over a domain  $\Delta \subset X(\mathbb R)$ defined by polynomial inequalities with coefficients in $\mathbb Q$; 

- the functional case: (functional) periods are integrals of differential forms $\omega$ on semialgebraic domains $\Delta$, depending algebraically on a parameter $z$ (the constant field being $\mathbb C$). 
These multivalued functions satisfy fuchsian linear differential equations (Picard-Fuchs). 
They have been revisited by J. Ayoub, who proved that {\it after performing an algebraic change of variable, $\Delta$ can be taken to be the unit polydisk}. 

 \smallskip Relations between periods have sparked speculations since Leibniz's time, \cf \cite{Wu}. In the functional case, the problem has been completely settled by Ayoub \cite{Ay1}\cite{Ay2} (\cf also \cite{BT}). In the arithmetic case, it remains very much open (conjectures of Grothendieck and Kontsevich-Zagier, \cf \cite[IX]{A1}\cite[XIII]{HW}\cite{KZ}).
    
  \subsection{G-functions and functional periods} Let us consider the situation where the functional periods are integrals on the fiber of a smooth pencil $X\to S$ defined a number field $K \subset \mathbb C$, $z$ being a local coordinate on $S$. Then the Picard-Fuchs differential operator which controls the variation of the periods is a G-operator, whence a relation between functional periods 
  and G-functions.  
    Namely, writing the operator in matrix form $d/dz - A(z)$, and assuming for simplicity that $z=0$ is not a singularity, there is a fundamental matrix $G(z)$ of solutions with $G(0) = I$ with G-function entries $G_{ij}(z)$; on the other hand, the period matrix $P(z)$ is a fundamental matrix of solutions, so that 
  $$ G(z) = P(z)P(0)^{-1}$$ (hence functional periods are linear combinations of G-functions with coefficients in the ring of arithmetic periods). 
  
  In particular, for any specialization $z\leadsto \zeta \in \bar K$ where $G$ converges, the G-values $G_{ij}(\zeta)$ are quotients of (arithmetic) periods. 
  One can be more precise, using Ayoub's reduction of the domain to a polydisk: {\it any period is a G-value} (Fres\'an, unpublished). The question of whether the converse holds is closely related to the Bombieri-Dwork conjecture.

  \subsection{Relations between periods, and relations between G-values} 
  \subsubsection{} If for some point $s\in S(\bar K)$, the fiber $X_s$ has more ``symmetries" than the general fiber \eg if it has more endomorphisms, in the usual or in the motivic sense - more precisely, if the motivic Galois of $X_s$ is smaller than the motivic Galois of the general fiber - , then this yields polynomial relations between the periods $P_{ij}(\zeta)$ of $X_s$, which do not come from polynomial relations between the functions $P_{ij}(z)$ (and according to Grothendieck's period conjecture, all such relations should come in this way). 
  
  If moreover there are sufficiently many such relations, elimination of the entries of $P(0)$ may yield polynomial relations between the G-values $G_{ij}(\zeta)$, which do not come from relations between the G-functions $G_{ij}(z)$. 
  
   \subsubsection{} A nice such setting arises when, instead of taking for $0$ a non-singularity as above, one considers certain columns of $G$ which are related to the notch $0$ of the monodromy-weight filtration (\ie if one looks at $n$-forms, the image of $(N/2\pi i)^n$, where $N$ is the logarithm of the unipotent part of the local monodromy at $0$). Let's call them {\it $M_0$-G-functions and $M_0$-periods} respectively.
   
 \smallskip  Fact: {\it $M_0$-G-functions are quotients of $M_0$-periods by $(2\pi i)^n$}. 
 
  \smallskip    A particularly interesting case, for $n=1$, is given by abelian pencils with multiplicative reduction at $z=0$. The simplest instance is the case of the Legendre elliptic pencil: a basis of $M_0$-periods is then given by $2\pi i F(1/2, 1/2, 1; z),  \, 2\pi i F(1/2, -1/2, 1; z)\,$.
   
   \smallskip    We refer to \cite[IX, Th. 2 and 1, \S 4]{A1} for precise definitions, statements and proofs\footnote{proofs are written in \cite{U} and in greater generality; see also \cite[\S 2]{P4}}.   
      
    \medskip  In this setting, if one can eliminate $2\pi i$ from relations between evaluations at $\zeta$ of $M_0$-periods, one gets relations between the corresponding G-values.

 \subsubsection{} Once one obtains polynomial relations between $G_{ij}(\zeta)$'s, one is in position to apply various levels of the diophantine theory of G-values, in the hope to get some diophantine constraints on $\zeta$. 
 
\smallskip The general theory only yields rather weak constraints \cf \cite[\S 12]{B}\cite[VII \S 4]{A1} (very recent progress came from using global properties of G-operators - arithmetic holonomy bounds \cite{CDT}). 
  
 But in very special instances, a method introduced in \cite{A4} provides optimal results; it applies to the G-functions $F(1/2, 1/2, 1; z),  F(1/2, -1/2, 1; z)$ and shows that the period matrix of an elliptic curve has transcendence degree at least 2; in the case of complex multiplication, this bound is sharp: it provides a {\it G-function-theoretic proof of Chudnovsky's theorem} (the period matrix of a CM elliptic curve has transcendence degree $2$).  
  
 \smallskip Bombieri's principle of global relations \cite[\S 11]{B} lies at an intermediate level: it gives the finiteness of the set of $\zeta$'s of bounded degree where there is an exceptional ``global relation" between G-values at $\zeta$. 
  
     \section{Bombieri's principle of global relations between G-values, and applications.} 
   
    \subsection{Bombieri's principle of global relations}  Let $\vec g(z) = (g_1, \ldots, g_m)$ be a sequence of G-functions satisfying a linear differential system $$(\frac{d}{dz} - A(z)) \vec g(z)  =0,$$ with $A_{ij}(z)\in K(z)$. Let $\zeta\in \bar K$.  
    
    A homogeneous polynomial relation $p( g_1(\zeta), \ldots, g_m(\zeta))=0$ of degree $\delta$ is called {\it global} if it holds $v$-adically for any place $v$ of $K(\zeta)$ such that $\vert \zeta\vert_v < \min (1, R_v(g))$. It is called {\it exceptional}\footnote{we prefer the terminology ``exceptional" to ``non-trivial": if $g$ is as before, so is $g_\zeta := (z-\zeta)g$ for any $\zeta\in K$; the relation $g_\zeta(\zeta)=0$ can hardly be called non-trivial} if it does not come by specialization from a homogeneous polynomial relation $q(g_1(z), \ldots, g_m(z))=0$ of degree $\delta$ with coefficients in $\bar K[z]$. 
    
  The principle of global relations 
   asserts that, {\it given $g(z)$ and $\delta$, the absolute logarithmic height of the set of $\zeta \in \bar K$ where an exceptional global relation of degree $\delta$ occurs between $g_1(\zeta), \ldots, g_m(\zeta)$ grows at most polynomially with $\delta$}. 
    
    In particular the subset consisting of such $\zeta$'s of bounded degree is finite. 
 
     \subsection{First applications}  \subsubsection{} The first applications of this principle in the context of Diophantine Geometry, specifically to relations between G-values coming from geometry (in the wake of 2.3.1 and 2.3.2 above) were given in \cite[X]{A1}\footnote{elaborating on the author's thesis}.     
    
  By combining the technique explained in 2.3.2 with the principle of global relations\footnote{in this situation $\delta$ is essentially proportional to $ [K(\zeta):K]$, and only archimedean places contribute to the global relation}, the following was proven in \cite[X, \S 1.3, 4.2]{A1}:
   {\it   let $S$ be an algebraic curve defined over $K$ lying on the moduli space $\sA_g$ of p.p.a.v. of odd dimension $g>1$, whose closure $\bar S$ in the Baily-Borel compactification meets the $0$-dimensional boundary component $\infty$ (\ie the corresponding abelian pencil has multiplicative reduction), and such that the abelian variety $A_\eta$ over the generic point of $S$ is absolutely simple. Then there are only finitely many {\rm special points}\footnote{\ie points corresponding to abelian varieties $A_s$ with complex mutiplication} $s\in S(\bar K)$ of bounded degree}. 
    
  \subsubsection{}   As explained in \cite[X \S 4.3]{A1}, it was tempting to think, by analogy with the Manin-Mumford conjecture (Raynaud's theorem), that the bunch of assumptions \{$g$ odd, $\bar S$ meeting $\infty$,  $A_\eta$ absolutely simple, and $s$ of bounded degree\} was irrelevant, provided $S$ is not a special curve. This was one of the sources of inspiration for the Andr\'e-Oort conjecture. 
     
     This line of thought also led to (a special case of) the Andr\'e-Pink conjecture \cite[X \S 4.5]{A1}\cite{Pi}.
    
     \subsubsection{} Let us summarize the technique (the {\it G-function method}): given a projective smooth fibration $X\to S$ and an exceptional point $s\in S(\bar K)$ for which the fiber $X_s$ has extra motivic symmetries, one gets period relations at $s$. In good situations (when one can eliminate periods of the fiber at the base point), one gets relations between G-values, and even global relations (for this, it suffices to construct separately polynomial relations between $v$-adic evaluations for each relevant $v$ and then to multiply them, if one can bound the number of such $v$'s in order to bound the degree $\delta$ - a crucial issue in this approach). 
     
     When the method applies, one thus concludes that {\it the set of such points $s$ of bounded degree, where there are extra motivic motivic symmetries in $X_s$, is finite}. 
     
  \subsubsection{}   Actually, A. Cadoret and A. Tamagawa {\it obtained this conclusion unconditionally} \cite{CT}\cite{C}, using the $\ell$-adic realization of motives instead of the period realization (at points where there are extra motivic symmetries, the $\ell$-adic arithmetic monodromy group is smaller than the generic one). 
  
     \smallskip      However, the G-function method, when applicable, gives more: namely, {\it a bound for the height of such exceptional points}. This turns out to be crucial in the new applications of this method.

    \subsection{Combination with the Pila-Zannier method, and new applications}(General references: \cite{Ch}\cite{Z})\label{s3} 
    \subsubsection{} A decisive innovation in the study of the Andr\'e-Oort and of the (more general) Zilber-Pink conjecture in the context of unlikely intersections has been the Pila-Zannier method (and before, the pioneering work \cite{BMZ}).  
   
   In this context, one deals with a parameter set $S$ (an algebraic $K$-variety) and a set $\Sigma\subset S(\bar K)$ of ``exceptional points", and one wishes to show that $\Sigma$ is not Zariski-dense. The hardest point in the Pila-Zannier method is to bound from below the size of the Galois orbit of points $\zeta \in \Sigma $ - which can often be reduced to showing that the absolute logarithmic height $ h(\zeta) $ grows at most polynomially with the degree $[K(\zeta): K]$.  
  
  When $S$ is a curve, the G-function method explained above, applied to the fibration parametrized by $S$, offers a possible approach. 
        
     \subsubsection{} However, the application of the principle of global relations requires to consider not only relations between complex G-values but also relations between $p$-adic G-values. 
     
     The $p$-adic counterpart of the setting of 2.3.2 ($M_0$-G-functions and $M_0$-periods) was established by D. Urbanik \cite[Th. 1.14]{U}: the $p$-adic functional periods are now understood in the sense $p$-adic Hodge theory, the needed integral structure on the relevant $M_0$-subspace on the $p$-adic (pro-)\'etale side comes from a ``graph-theoretic" homology group\footnote{actually Urbanik proceeds a little differently in obtaining and applying his result}. Denoting by $t_p\in \mathbb Z_p(1) \subset B_{dR}$ a $p$-adic analog of $2\pi i$, one gets the $p$-adic analog of 2.3.2:
     
  \smallskip Fact: {\it $M_0$-G-functions are quotients of $M_0$-periods by $t_p^n$}.   
     
   \smallskip     The special case of an abelian pencil with multiplicative reduction was already considered in \cite[\S 5.2.2, 5.2.3]{A5}\footnote{via \cite{A2}, which relied however upon unpublished results by Raynaud}; other cases of abelian pencils were treated independently by C. Daw and M. Orr \cite{DO1}. 
     
   \subsubsection{}  A number of applications of this method to various special cases of the Zilber-Pink conjecture in the same spirit appeared around the same time. Let us mention a sample.
   
 To begin with, this improvement of the G-function method with finite places coming into play allows Daw and Orr to remove, in 3.2.1 above, the assumptions \{$g$ odd and $A_\eta$ absolutely simple\} when $S$ is Hodge-generic\footnote{this means that the generic Mumford-Tate group (or motivic Galois group) is maximal, \ie the group of symplectic similitudes $CSp_{2g}$. Equivalently, here: the monodromy group is Zariski-dense in $Sp_{2g}$}. Via the Pila-Zannier machinery, they get a {\it new proof of the Andr\'e-Oort conjecture for Hodge generic curves $S$ in $\sA_g$ ($g>1$) whose closure $\bar S$ meets $\infty$}\footnote{this condition is somewhat relaxed in Papas' work \cite{P2}; slightly earlier,  G. Binyamini and D. Masser got the Andr\'e-Oort conjecture for non-compact curves in Hilbert modular varieties by a similar technique \cite{BM}}.
   
 They use the same method to prove the {\it Zilber-Pink conjecture for Hodge generic curves $S$ in $\sA_2$ such that $\bar S$ meets $\infty$: the intersection of $S$ with the union of special curves in $\sA_2$ is finite}\footnote{we take this opportunity to mention that: $i)$ Problem 2 in \cite[X \S 4.4]{A1} has a negative answer; the correct statement of this kind is the Zilber-Pink conjecture; $ii)$ in Problem 3 of \cite[X \S 4.5]{A1}, ``Shimura subvariety" should be replaced by ``weakly special curve"; $iii)$ \cite[X app.]{A1} is flawed (it is not true that $\bar\tau_v -\tau_v $ takes only two values), but the claim that one can prove the transcendence of $\pi$ by G-function methods is substantiated in \cite{A4}} \cite{DO1} (\cf also \cite[\S 1.5]{P2}). 
   
  Further, they prove the {\it Zilber-Pink conjecture for Hodge-generic\footnote{in this case, this is equivalent to ``$S$ is not contained in a special subvariety"} curves $S$ in some power $Y(1)^n$ of the modular curve, such that $\bar S$ meets $\infty$ (\ie the corresponding abelian pencil has multiplicative reduction): the intersection of $S$ with the union of special subvarieties of codimension $>1$ is finite} \cite{DO2}\footnote{here the condition that $\bar S$ meets $\infty$ replaces the ``asymmetry" condition in earlier work by P. Habegger and J. Pila \cite{HP}, and is further somewhat relaxed by Papas \cite{P3}}.
   
   \smallskip Urbanik considers Hodge-generic curves $S$ on the moduli space $\sM_g$ of curves of genus $g>1$, whose closure $\bar S$ in the Deligne-Mumford compactification meets the boundary in the locus of stable curves with components of genera $g_i$ and $\sum g_i \leq g-2$. He proves that {\it the intersection of $S$ with the set $\Sigma \subset \sM_g$ of points corresponding to curves whose jacobian has an isogeny factor with complex multiplication is finite} \cite[Th. 1.1]{U}.  
   
 \smallskip Pursuing this approach, G. Papas gets results of the kind but for different conditions at the boundary, and for sets $\Sigma$ defined by different conditions about the endomorphisms of the corresponding jacobians \cite{P4}.  
  
  Both Urbanik \cite{U} and Papas \cite{P4}   obtain much more general statements concerning atypical intersections \`a la Klingler \cite{Kl} in the context of variations of Hodge structures. 
   
       \subsubsection{} In all aforementioned applications of the G-function method to the Zilber-Pink conjecture, there is a condition about $\bar S \setminus S$ ($z$ is a local parameter around a point of $\bar S \setminus S$). This is mainly due to the following reason. 
       
       In the G-function method, one wishes to show that exceptional points $s\in \Sigma$ give rise to relations between relevant G-functions evaluated at $\zeta = z(s)$  for (essentially) every place $v$ where the evaluation makes sense. 
     The method, explained in 2.3.1 for an archimedean place $v$, is based on first establishing period relations and then passing from periods to G-values. For a non-archimedean place $v$, a workable notion of $v$-adic periods, relating $v$-adically two arithmetic structures, is constructed in the aforementioned cases only in the setting of {\it bad reduction} at $v$, typically in the setting of 2.3.2 (when $\zeta$ is $v$-adically small, \ie $s$ closed to $\bar S \setminus S$).  
     
 Yet there are cases of good reduction where one can indeed define a workable notion of $p$-adic periods, relate them to $p$-adic G-values, and get in this way polynomial relations between $p$-adic G-values: this is the case for an abelian pencil $A\to S$ in the $p$-adic locus of {\it supersingular reduction}. We refer to \cite[\S 5.3]{A5} and \cite{A3} for a detailed discussion of the case of an elliptic pencil with CM fibers $A_0, A_s$ (in this case, a formula due to Gross-Zagier-Dorman-Lauter-Viray \cite{LV} gives some grasp on the number of $v$'s such that  $\zeta= z(s)$ is $v$-adically small); see also \cite{AF}. 
  
    Could this enlargement of the scope of the G-function method in the $p$-adic direction lead to new cases of Zilber-Pink?   
    
          \subsubsection{} In the previous results, the G-function method plays only one part in the Pila-Zannier strategy for unlikely intersection problems. The Pila-Wilkie counting theorem provides the main tool, leading to powerful transcendence theorems about functional periods and beyond. 
     
           One may wonder whether this powerful tool could also lead to transcendence theorems about numbers, \eg (arithmetic) periods\footnote{a test bench could be Schneider's problem: to prove, {\it without using elliptic functions}, that for a non-CM elliptic curve defined over $\bar \Q$, the corresponding period quotient $\tau\in \mathfrak H$ is transcendental - viewed as a bialgebraicity problem involving the $j$-function (but not elliptic functions).}.  
        
        \bigskip\noindent
        {\smit Acknowledgements: I am grateful to J. Fres\'an for a list of useful comments.}

\end{sloppypar} 

 \end{document}